\documentclass{amsart}
\usepackage{amssymb}
\usepackage{xfrac} 
\usepackage{comment}
\usepackage{mathtools}
\usepackage{tikz-cd}
\usepackage{bm}
\usepackage{enumerate}
\usepackage{enumitem}

\usepackage{amsmath,amssymb,amsthm,amsfonts,mathrsfs}

\usepackage[backend=biber,style=alphabetic,maxnames=10,minalphanames=3]{biblatex}
\appto{\bibsetup}{\raggedright}

\usepackage{caption}
\usepackage{subcaption}

\usepackage[colorlinks,cite color=blue,pdftex,linkcolor = purple]{hyperref}
\usepackage{cleveref}

\newtheorem{theorem}{Theorem}[section] 
\newtheorem{lemm}[theorem]{Lemma}
\newtheorem{proposition}[theorem]{Proposition}

\newtheorem{corollary}[theorem]{Corollary}

\theoremstyle{definition}

\newtheorem{definition}[theorem]{Definition}
\newtheorem{remark}[theorem]{Remark}

\DeclareMathOperator{\rk}{rk}

\newcommand{\cl}{\operatorname{cl}}
\newcommand{\trdeg}{\operatorname{trdeg}}

\newcommand{\RR}{\mathbb{R}}
\newcommand{\ZZ}{\mathbb{Z}}
\newcommand{\FF}{\mathbb{F}}
\newcommand{\QQ}{\mathbb{Q}}
\newcommand{\Ga}{\mathbb{G}_{\mathrm a}}

\newcommand{\NL}{\mathrm{NL}}
\newcommand{\Stab}{\operatorname{Stab}}
\newcommand{\Spec}{\operatorname{Spec}}
\newcommand{\Hom}{\operatorname{Hom}}

\begin{document}
\title{Algebraic matroids with non-algebraic duals exist in any positive characteristic}

\author{Lukas Grund}
\address{Lukas Grund \\ Friedrich-Schiller-Universit\"at Jena, Fakult\"at f\"ur Mathematik und Informatik, Institut fur Mathematik, Inselplatz 5, 07743 Jena, Germany}
\email{lukas.grund@uni-jena.de}

\begin{abstract}
We show that for every field $K$ of positive characteristic there exists an algebraic matroid $M$ over $K$ such that the dual matroid $M^*$ is not algebraic over any field. Our conceptual proof relies on the theory of matroids over one-dimensional algebraic groups. As an application, we deduce that the classes of supports of realizable volume polynomials and of realizable covolume polynomials differ in positive characteristic.
\end{abstract}

\maketitle

\section{Introduction}\label{sec:introduction}

Matroids were introduced by Whitney \cite{whitney} to put the similar properties of subsets of linearly independent vectors and trees in a graph into an abstract combinatorial framework. It was already noted earlier by van der Waerden that algebraically independent subsets in a field extension satisfy the same axioms~\cite{van1931moderne}. Concretely, let $M=(E,r)$ be a \emph{matroid} with a finite ground set $E$ and a \emph{rank function} $r$, see e.g.\ \cite{oxley}.
\begin{definition}
	A matroid $M=(E,r)$ is \emph{algebraic} over a field $K$ if there exists a field extension $K \subseteq L$ and $e_i\in L$ for $i\in E$ such that for all $I\subseteq E$ it holds that
	\[
		r(I)=\trdeg_{K}K(e_i : i\in I).
	\]
\end{definition}
Here, $\trdeg_K$ denotes the transcendence degree over $K$ and $K(e_i : i\in I)$ is the smallest subfield of $L$ containing $K$ and $\{e_i: i\in I\}$.
Hence, $I\subseteq E$ is \emph{independent} in $M$, that is $r(I)=|I|$, if and only if the set $\{e_i: i\in I\}$ is \emph{algebraically independent} over $K$.

If a matroid is representable over some field $K$, it is also algebraic over $K$.
Lindstr\"om showed in \cite{lindstrom1989matroids} that being algebraic over $K$ depends only on the characteristic of $K$. In characteristic $0$, every algebraic matroid over $K$ is representable over a field extension of $K$ \cite{Ing71}. In positive characteristic however, Ingleton described a matroid on $11$ elements of rank $3$ which  is algebraic over $\FF_2$  and has both the Fano and non-Fano matroids as restrictions; hence it is not representable over any field.
It was only shown in 1975 that not all matroids are algebraic: Ingleton and Main proved that the V\'amos matroid is not algebraic over any field~\cite{IM75}.

In general, algebraic matroids are far less understood than linearly representable matroids. It was recently shown by Boege and Yashfe \cite{boege2026recognition} that the problem of recognizing whether a given matroid is algebraic is undecidable.
Algebraic matroids respect several prominent matroid operations.
The class of algebraic matroids over a field $K$ is closed under taking minors, truncations, principal extensions, matroid unions, and matroid intersections~\cite{GHMSW2025preservingvolume,matuvs2018classes}.

Recall that the \emph{dual} $M^*$ of a matroid $M$ is the matroid on the same ground set with bases the complements of bases of $M$. This notion is for instance motivated by the \emph{orthogonal complement} of a vector space and by \emph{plane graph duality}, since $M(G^*)=M^*(G)$ for a planar graph $G$, where $M(G)$ is the graphic matroid of a graph $G$. Whether the dual of an algebraic matroid is again algebraic in general was recently answered negatively by Larson and Le~\cite{larson2026duals}. However, their counterexample only works over a field of characteristic $2$ and the proof relies on a computer calculation. Inspired by their counterexample, we construct counterexamples over every field of positive characteristic and give a conceptual proof. For this we use generalizations of the \emph{Fano} and \emph{non-Fano} matroids, called the \emph{Lazarson matroid} $L_p$ and the \emph{non-Lazarson} matroid $\NL_p$ over $\overline{\FF}_p$. Let $u_1,\dots,u_{p+1},x_1,\dots,x_{p+1}$ be algebraically independent over $\overline{\FF}_p$. Then $L_p$ is represented by
\[x_1,\dots,x_{p+1},y_1,\dots,y_{p+1},z\]
where $z = \sum_{i=1}^{p+1} x_i$ and $y_i = z - x_i$.
$\NL_p$ is represented by
\[u_1,\dots,u_{p+1},v_1,\dots,v_{p+1},w\]
where $w = \prod_{i=1}^{p+1} u_i$ and $v_i = \frac{w}{u_i}$. 
We show that every algebraic realization of the dual matroids $L_p^*$ and $\NL_p^*$ is algebraically equivalent to a realization over a one-dimensional algebraic group: over $\mathbb{G}_a$ for $L_p^*$, and over $\mathbb{G}_m$ or an elliptic curve for $\NL_p^*$. It follows that every algebraic matroid $N$ obtained by gluing $L_p^*$ and $\NL_p^*$ along a common basis is the intersection of $L_p^*$ and $\NL_p^*$, each extended by coloops. Dually, $N^*$ is the union of $L_p$ and $\NL_p$, each extended by loops. We then construct an algebraic matroid $M_p$ that has both $L_p$ and $\NL_p$ as contractions but is not such a union. Therefore the dual of $M_p$ is not algebraic over any field.

\begin{theorem}[Theorem \ref{thm:counterexample}]\label{thm:main}
	Over any field $K$ of \emph{positive characteristic} there is an \emph{algebraic matroid} $M$ over $K$ such that the \emph{dual matroid} $M^*$ is \emph{not algebraic} over any field.
\end{theorem}

Aluffi~\cite{aluffi2024lorentzian} introduced the notion of \emph{covolume polynomials} which were further developed in~\cite{GHMSW2025preservingvolume}. Matroids whose set of bases forms the \emph{support} of some \emph{multiaffine realizable covolume polynomial} over a field $K$ are exactly those matroids $M$ whose \emph{dual} $M^*$ is an \emph{algebraic matroid} over $K$. To any covolume polynomial over $K$ we can associate a \emph{volume polynomial} over $K$ by generalizing matroid duality.

\begin{theorem}[Theorem \ref{thm:volume}]
    Let $K$ be a field of \emph{positive characteristic}. There exists a \emph{realizable covolume polynomial} $g$ over $K$ such that the support of $g$ is not the support of any \emph{realizable volume polynomial}.
\end{theorem}

\subsection*{AI disclosure} LLMs were used for proofreading and editing. Claude Opus 5 and 5.5 assisted us in finding the statement and the proof of Theorem \ref{thm:matroid-intersection}.

\subsection*{Acknowledgements}
The author thanks Giles Gardam, June Huh and Lukas Kühne for the many insightful discussions on algebraic matroids. 
The author was supported by the SPP 2458 "Combinatorial Synergies".

\section{Matroids over one-dimensional groups}
Throughout the whole article let $K$ be an algebraically closed field of positive characteristic $p$ and $L$ an algebraically closed field extending $K$ of infinite transcendence degree over $K$.

We will use the theory of matroids over a \emph{connected one-dimensional algebraic group} $G$ over $K$ first described by Evans and Hrushovski~\cite{evans1991projective} and further developed in~\cite{BCD}. An \emph{algebraic group} is said to be \emph{connected} if it is \emph{irreducible}. All algebraic groups appearing in this article will be connected. Let $\mathbb{E}$ denote the \emph{endomorphism ring} of $G$. This is a possibly non-commutative ring without zero divisors and an Ore domain. Let $Q$ be the corresponding \emph{division ring} of $\mathbb{E}$. In all the commutative algebraic groups appearing we denote the addition with $\oplus$. For more on algebraic groups see~\cite{milne2017algebraic}. There are the following equivalent definitions for a \emph{matroid over a one-dimensional group}:
\begin{definition}
\label{def:matroid-over-group}
    \begin{enumerate}[label=(\alph*), ref=(\alph*)]
        \item \label{def:matroid-over-group:1} Given a closed connected subgroup $X$ of $G^n$ of dimension $d$. We write $M(X)$ for the matroid whose independent sets are all the sets $I \subseteq [n]$ such that the projection of $X$ to $G^I$ is all of $G^I$.
        \item \label{def:matroid-over-group:2} Given a homomorphism of algebraic groups $\psi: G^d \to G^n$. Then $\psi$ is encoded by a matrix $\psi \in \mathbb{E}^{n \times d}$. Let $M(\psi)$ be the matroid defined by the (left) row space of $\psi$ over $Q$.
        \item \label{def:matroid-over-group:3} Given a generic $L$-point $p:\Spec L \to G^d$. The image $\psi(p)$ gives $n$ points $q_1,\dots,q_n \in G(L)$. We denote with $M(q_1,\dots,q_n)$ the algebraic matroid realized by those points. 
    \end{enumerate}
\end{definition}

I.e. we have $M(X) = M(\psi) = M(q_1,\dots,q_n)$. The equivalence of (a) and (b) follows from \cite[Theorem 1.1]{BCD} and the equivalence of (a) and (c) is easy to see.
For any $p \in G(L)$ we write $\cl(p)$ for the \emph{algebraic closure} of $K(p)$. Here $G(L) = \Hom_K(\Spec L, G)$ and $p$ corresponds to a point $g\in G$ together with an embedding $\kappa(g)\hookrightarrow L$; let $K(p)\subseteq L$ be its image.

The one-dimensional algebraic groups over $K$ are completely classified. Such a group $G$ is either the \emph{additive group} $\mathbb{G}_a$ over $K$, the \emph{multiplicative group} $\mathbb{G}_m$ over $K$ or an \emph{elliptic curve} $E$ over $K$. See~\cite[Theorem III.10.9]{borel2012linear} for a classification of affine connected one-dimensional algebraic groups. Thus the possible \emph{endomorphism rings} are understood quite well. The \emph{endomorphism ring} of $\mathbb{G}_m$ or an elliptic curve $E$ has characteristic $0$ while the \emph{endomorphism ring} of $\mathbb{G}_a$ has the same characteristic as $K$.
Matroids over one-dimensional groups share several properties with \emph{linear matroids} over fields. For example, they can always be extended \emph{modularly}.

\begin{lemm}
\label{lemm-modular-intersection}
    Let $G$ be a one-dimensional group over $K$ with endomorphism ring $\mathbb{E}$ and division ring $Q$. For any $p \in G(L)^d$ \emph{generic} and any $v \in \mathbb{E}^{d}$ write $v(p) = \bigoplus_j v^j(p_j) \in G(L)$. If we are given $u_1,\dots,u_n,v_1,\dots,v_m \in \mathbb{E}^{d}$ such that the (left) $Q$-spans $U = \langle u_1,\dots,u_n \rangle$ and $V = \langle v_1,\dots,v_m \rangle$ have dimension $n$ and $m$ and $U + V$ has dimension $n+m-k$, then $U \cap V$ is spanned by some $w_1, \dots, w_k \in \mathbb{E}^{d}$ and the points $q_i = w_i(p)$ are \emph{generic} and \emph{independent} with
    \[\cl(u_1(p),\dots,u_n(p)) \cap \cl(v_1(p),\dots,v_m(p)) = \cl(q_1,\dots,q_k).\]
\end{lemm}
\begin{proof}
    We get the $w_i$ through the dimension formula for left vector spaces over $Q$. As matroids over $Q$ are equivalent to algebraic matroids over $G$ the statement follows.
\end{proof}
We state some properties of algebraic groups which we will need later.
\begin{lemm}
\label{lemm:matroids-over-groups}
    Let $G$ be a \emph{connected one-dimensional algebraic group} over $K$. 
    \begin{enumerate}[label=(\alph*), ref=\alph*]
        \item \label{lemm:matroids-over-groups:1} The \emph{action} of any $e \in \mathbb{E} \setminus \{0\}$ \emph{preserves the closure} of a point $p \in G(L)$. We have
        \[\cl(p) = \cl(e(p)).\]
        \item \label{lemm:matroids-over-groups:2} Let $H$ be a commutative connected algebraic group. For subvarieties $V,C$ we write $V \oplus C \coloneqq \{v \oplus c \mid v \in V, c \in C\}$. Let $V$, $C$ be closed and irreducible and $W = \overline{V \oplus C}$. If $\dim W = \dim V$, then $C \ominus C \subseteq \Stab (W)$ and $W = V \oplus C$.
       \item \label{lemm:matroids-over-groups:3} If $G$ is either an \emph{elliptic curve} or $\mathbb{G}_m$, then every \emph{connected closed subgroup} of $\mathbb{G}_a^n \times G^m$ is a product $U \times V$ where $U$ is a subgroup of $\mathbb{G}_a^n$ and $V$ is a subgroup of $G^m$.
    \end{enumerate}
\end{lemm}
\begin{proof}
    (a) Consider the composition of morphisms $\Spec L \xrightarrow{p} G \xrightarrow{e} G$. We get $L \hookleftarrow K(p) \hookleftarrow K(e(p))$ and as $e$ defines an isogeny the extension $K(p)/K(e(p))$ is finite. Hence $\cl(p) = \cl(e(p))$. \\ (b) $V \times C$ is irreducible and $V \oplus C$ is the image under the morphism $\oplus: H \times H \to H$. Hence $W$ is also irreducible.
    For any $c \in C(K)$ we have $V \oplus c \subset V \oplus C \subseteq W$ and $\dim(V \oplus c) = \dim V = \dim W$. With both varieties being closed and irreducible it follows $V \oplus c = W$ for every $c \in C(K)$. For $c,c' \in C(K)$ we now get
    \[(c \ominus c') \oplus W = (c \ominus c') \oplus (c' \oplus V) = c \oplus V =  W\]
    and with that $c \ominus c' \in \Stab (W)$. It follows $C \ominus C \subseteq \Stab (W)$ and $W = V \oplus C$. \\
    (c) The statement is standard. Let $W$ be a closed connected subgroup of $\mathbb{G}_a^n \times G^m$. For $G = \mathbb{G}_m$ it follows from the unique decomposition of $W$ into a product of a \emph{unipotent} subgroup and a subgroup of \emph{multiplicative type}, see~\cite[Theorem 16.13]{milne2017algebraic}. For $G = E$ an elliptic curve it follows from the \emph{Rosenlicht decomposition}, see~\cite[Theorem 6.2]{milne2013proof} and \cite{rosenlicht1956some}, that $W = (W_{\text{aff}} \times W_{\text{ant}}) / (W_{\text{aff}} \cap W_{\text{ant}})$ in an affine and anti-affine group.
    Here $W_{\text{aff}} \subset \mathbb{G}_a^n \times 0$ and $W_{\text{ant}} \subseteq 0 \times E^m$. As these intersect trivially, we get the wanted result with $U = W_{\text{aff}}$ and $V = W_{\text{ant}}$.
\end{proof}

To study the algebraic independence of elements in a field extension geometrically, we define the following:
\begin{definition}
    Let $V$ be a variety over $K$ and $p\colon \Spec L\to V$ be an $L$-point. We call the irreducible closed subvariety 
    \[Z(p)\coloneqq \overline{\{p\}}\subseteq V \]
    the \emph{locus} of $p$.
\end{definition}

\begin{lemm}
\label{lemm:loci}
    Let $V,W$ be varieties over $K$ and $p \in V(L)$ an $L$-point.
    \begin{enumerate}[label=\textup{(\alph*)}, ref=\textup{\alph*}]
        \item \label{lemm:loci:1} $\dim Z(p)=\trdeg_KK(p)$.
        \item \label{lemm:loci:2} If $\varphi:V\to W$ is a morphism over $K$, then $Z(\varphi(p))=\overline{\varphi(Z(p))}$.
        \item \label{lemm:loci:3} If we are given several points $p_i\in V_i(L)$ in varieties $V_1,\dots,V_n$ over K and \\ $\trdeg_K K(p_1,\dots,p_n)=\sum_i\trdeg_K K(p_i)$, then the locus splits as a product of loci, \\ i.e. $Z(p_1,\dots,p_n)=Z(p_1)\times\dots\times Z(p_n)$.
    \end{enumerate}
\end{lemm}
\begin{proof}
    (a) The dimension of $Z(p)$ is just the transcendence degree of the residue field of its generic point. This is $\trdeg_KK(p)$.
    (b) $\overline{\varphi(Z(p))}$ is closed and contains $\varphi(p)$. Hence $Z(\varphi(p)) \subseteq \overline{\varphi(Z(p))}$. The other inclusion follows by continuity. 
    (c) We set $p = (p_1,\dots,p_n)$. With (b) we get $Z(p) \subseteq \prod_i Z(p_i)$. Both sides are now irreducible and of the same dimension, hence we get an equality.
\end{proof}

In Section \ref{sec:alg-reps} we glue together different realizations of the matroid $M(K_4)$ of the complete graph on four vertices. For this we will need the following proposition which is a variant of \cite[Theorem 2.2.2]{evans1991projective}.
\begin{proposition}
\label{lemm-gauge}
    Let $G_1,G_2$ be one-dimensional algebraic groups over $K$. If for $a_1,b_1 \in G_1(L)$ and $a_2,b_2 \in G_2(L)$ we have
    \begin{enumerate}
        \item $a_1$ and $b_1$ are \emph{algebraically independent},
        \item $\cl(a_1) = \cl(a_2)$ and $\cl(b_1) = \cl(b_2)$,
        \item $\cl(a_1 \oplus  b_1) = \cl(a_2 \oplus b_2)$.
    \end{enumerate}
    Then there is a closed connected one-dimensional subgroup $H_0 \subset G_1 \times G_2$ and points $p_a,p_b \in H_0(L)$ such that there are constants $x_a,x_b \in (G_1 \times G_2)(K)$ with
    \[a_1 =  \pi_1(p_a  \oplus x_a), \ a_2 = \pi_2(p_a \oplus x_a), \ b_1 =  \pi_1(p_b \oplus x_b),\ b_2 = \pi_2(p_b \oplus x_b).\]
    Here $\pi_1,\pi_2$ are just the coordinate projections of $G_1 \times G_2$. In addition $G_1$ and $G_2$ are \emph{isogenous}.
\end{proposition}
\begin{proof}
    Set $H = G_1 \times G_2$. Set $a = (a_1,a_2)$, $b = (b_1,b_2)$, $C_a = Z(a)$ and $C_b = Z(b)$, then $C_a$ and $C_b$ are irreducible curves in $H$ which dominate both factors with Lemma \ref{lemm:loci}\ref{lemm:loci:2}. Using Lemma \ref{lemm:loci}\ref{lemm:loci:3} follows $Z((a,b)) = C_a \times C_b$ in $H \times H$; hence $Z(a \oplus b) = \overline{(C_a \oplus C_b)}$ with Lemma \ref{lemm:loci}\ref{lemm:loci:2}. Furthermore $(1)$ and $(3)$ give us $\trdeg K(a \oplus b) =1$.
    By Lemma \ref{lemm:matroids-over-groups}\ref{lemm:matroids-over-groups:2} follows $C_b \ominus C_b \subseteq \Stab (Z(a \oplus b))$ and $C_a \ominus C_a \subseteq \Stab (Z(a \oplus b))$. Set $H_0 = \Stab (Z(a \oplus b))^\circ$. By translating with some points $x_a \in C_a(K),x_b \in C_b(K)$ we get by using that $C_a$ and $C_b$ are irreducible 
    \[C_a \subseteq x_a \oplus H_0, \ C_b \subseteq x_b \oplus H_0.\]
    With $x_a \oplus x_b  \oplus H_0 \subseteq Z(a \oplus b)$ follows $x_a \oplus x_b  \oplus H_0 = Z(a \oplus b)$ and $\dim H_0 = \dim Z(a \oplus b) =1$. The image of an algebraic group under a homomorphism is closed and $C_a$ dominates both factors; hence $\pi_1(H_0) = G_1$ and $\pi_2(H_0) = G_2$. The kernels are zero-dimensional and hence finite, therefore $G_1$ and $G_2$ are isogenous through $H_0$. The points are now $p_a = a \ominus x_a$ and $p_b = b \ominus x_b$. 
\end{proof}

\section{Intersecting matroids over division rings}
An important operation on matroids in this article will be the \emph{intersection of matroids}. The operation is defined in a purely combinatorial way. For two linear matroids over a common field  the \emph{intersection} is still linear over a sufficiently large extension of that field and for algebraic matroids the intersection is algebraic over the same field. Two matroids algebraic over the same field $K$ can however be linear over fields or division rings of different characteristic. This is the situation we need to understand. The dual notion of intersecting matroids is the \emph{union of matroids}. The operations are defined in the following way:
\begin{definition}
    Let $N_1$ and $N_2$ be matroids on a common ground set $E$.
    \begin{enumerate}
        \item The union $N_1 \vee N_2$ is the matroid on $E$ whose \emph{independent} sets are the subsets of $E$ that
        can be expressed as a \emph{union} $I_1 \cup I_2$, where $I_1$ is an \emph{independent} set of $N_1$ and $I_2$ is an
        \emph{independent} set of $N_2$.
        \item The \emph{intersection} $N_1 \wedge N_2$ is the matroid on $E$ whose \emph{spanning} sets are the subsets of $E$ that can be expressed as an \emph{intersection} $S_1 \cap S_2$, where $S_1$ is a \emph{spanning} set of $N_1$ and $S_2$ is a \emph{spanning} set of $N_2$.
    \end{enumerate}
\end{definition}

The two constructions are related through matroid duality by $N_1 \wedge N_2 = (N_1^* \vee N_2^*)^*$.
We say that $N_1 \vee N_2$ has the \emph{expected rank} if the set of bases of $N_1 \vee N_2$ is
\[\{B_1 \cup B_2 \mid B_1 \in J_1, B_2 \in J_2, \text{ and } B_1 \cap B_2 = \emptyset\}, \]
where $J_1$ is the set of bases of $N_1$ and $J_2$ is the set of bases of $N_2$. Dually, we say that $N_1 \wedge N_2$ has the \emph{expected rank} if the set of bases of $N_1 \wedge N_2$ is
\[\{B_1 \cap B_2 \mid B_1 \in J_1, B_2 \in J_2, \text{ and } B_1 \cup B_2 = E\}.\]

Now let $M_1$ be a matroid of rank $d$ with ground set $E_1$ over a division ring $Q_1$ and $M_2$ a matroid of the same rank $d$ with ground set $E_2$ over a division ring $Q_2$. We assume that the ground sets intersect in the standard unit vectors $E_1 \cap E_2 = T =\{e_1,\dots,e_d\}$. Let $M_1'$ be $M_1$ with the elements of $E_2 \setminus T$ added as coloops, and define $M_2'$ in the same way for $M_2$. We set $N = M_1' \wedge M_2'$ with ground set $E = E_1 \cup E_2$. Set $A = E_1 \setminus T$ and $B = E_2 \setminus T$.
The dual of $N$ has rank $|A| + |B|$. It follows that $N$ has rank $d$ and therefore $N$ has expected rank. This means that the set of bases of $N$ is
\[\{B_1 \cap B_2 \mid B_1 \text{ is basis of } M_1', B_2 \text{ is basis of } M_2' \text{ and } B_1 \cup B_2 = E\}.\]

To characterize the bases of the intersection of two matroids over possibly different division rings we will use \emph{Rado's theorem}, see~\cite[Theorem 11.2.2]{oxley} and \cite{Rado}.
\begin{theorem}[Rado's theorem]
    Given a matroid $M$ with ground set $E$ and rank function $r$. Given further a index set $I$ and a family $ \mathcal{A} = (A_i)_{i\in I}$ of subsets $A_i \subseteq E$. Then the following are equivalent:
    \begin{enumerate}
        \item The family $\mathcal{A}$ has a \emph{transversal} which is an \emph{independent} set of $M$.
        \item For every subset $J \subseteq I$ we have $r(A(J)) \geq |J|$.
    \end{enumerate}
    Here we set $A(J) = \bigcup_{j\in J} A_j$.
\end{theorem}
Over a division ring $Q$, the left row rank of a matrix equals its right column rank; we denote it by $\rk_Q$.

In the setting above, let $M_1$ be the row matroid of a matrix $\alpha$ over $Q_1$ with rows indexed by $E_1$ and columns indexed by $T$, and let $M_2$ be the row matroid of a matrix $\beta$ over $Q_2$ with rows indexed by $E_2$ and columns indexed by $T$, such that the rows indexed by $T$ are the standard unit vectors in both matrices. Let $S\subseteq E$ with $|S|=d$, and for $I\subseteq T$ let $\alpha_{S,I}$ be the submatrix of $\alpha$ with rows $S\cap E_1$ and columns $I$, and $\beta_{S,I}$ the submatrix of $\beta$ with rows $S\cap B$ and columns $I$.
\begin{theorem}
\label{thm:matroid-intersection}
     $S$ is a basis of $N$ if and only if
    \[\rk_{Q_1} \alpha_{S,I} + \rk_{Q_2} \beta_{S,I} \geq |I| \quad\text{for all } I \subseteq T.\]
\end{theorem}
\begin{proof}
    Write $S_T=S\cap T$, $S_A=S\cap A$, $S_B=S\cap B$, $a=|S\cap E_1|=|S_T|+|S_A|$ and $b=|S_B|$, so that $a+b=d$. We show that both conditions are equivalent to
    \begin{enumerate}
        \item[$(\ast)$] there is a partition $T=I_1\sqcup I_2$ such that $\alpha_{S,I_1}$ and $\beta_{S,I_2}$ are square and invertible.
    \end{enumerate}

    \emph{The inequality is equivalent to $(\ast)$.} Let $C_1$ be the column matroid of $\alpha_{S,T}$ on $T\times\{1\}$ and $C_2$ the column matroid of $\beta_{S,T}$ on $T\times\{2\}$, and let $C=C_1\oplus C_2$. For $t\in T$ set $A_t=\{(t,1),(t,2)\}$. For $I\subseteq T$, the rank of $\bigcup_{t\in I}A_t$ in $C$ is $\rk\alpha_{S,I}+\rk\beta_{S,I}$. By Rado's theorem, the inequality holds if and only if there is a partition $T=I_1\sqcup I_2$ such that the columns $I_1$ of $\alpha_{S,T}$ and the columns $I_2$ of $\beta_{S,T}$ are independent. Then $|I_1|\le a$ and $|I_2|\le b$, and $|I_1|+|I_2|=d=a+b$ forces equality, so $\alpha_{S,I_1}$ and $\beta_{S,I_2}$ are square and invertible. The converse is clear.

    \emph{$(\ast)$ implies that $S$ is a basis of $N$.} If $t\in S_T\setminus I_1$, then the row $e_t$ of $\alpha_{S,I_1}$ vanishes, so $(\ast)$ gives $S_T\subseteq I_1$. Set $B_1'=S_T\cup S_A\cup I_2$, a set of size $a+b=d$. The rows of $\alpha$ indexed by $I_2$ are unit vectors, so eliminating them together with the columns $I_2$ shows that the rows $B_1'$ of $\alpha$ are independent if and only if $\alpha_{S,I_1}$ is invertible. Hence $B_1'$ is a basis of $M_1$. Similarly, $B_2'=S_B\cup I_1$ is a basis of $M_2$. Then $B_1=B_1'\cup B$ is a basis of $M_1'$ and $B_2=B_2'\cup A$ is a basis of $M_2'$, with $B_1\cap B_2=S$ and $B_1\cup B_2=E$. Hence $S$ is a basis of $N$.

    \emph{If $S$ is a basis of $N$, then $(\ast)$ holds.} Let $B_1$ and $B_2$ be bases of $M_1'$ and $M_2'$ with $B_1\cap B_2=S$ and $B_1\cup B_2=E$. Since $B$ consists of coloops of $M_1'$ and $A$ of coloops of $M_2'$, we have $B_1=B_1'\sqcup B$ and $B_2=B_2'\sqcup A$ for bases $B_1'$ of $M_1$ and $B_2'$ of $M_2$. Set $I_1=T\cap B_2'$ and $I_2=T\setminus B_2'$. Since $B\subseteq B_1$, we get $B_2'\cap B=S_B$, so $B_2'=S_B\cup I_1$. Since $A\subseteq B_2$ and $B_1\cup B_2=E$, we get $B_1'\cap A=S_A$ and $B_1'\cap T=S_T\cup I_2$, where $S_T\subseteq I_1$. Hence $B_1'=S_T\cup S_A\cup I_2$, and $d=|B_1'|=a+|I_2|$ gives $|I_2|=b$ and $|I_1|=a$. The same block elimination as above shows that $\alpha_{S,I_1}$ and $\beta_{S,I_2}$ are invertible.
\end{proof}

\begin{remark}
    We note that Theorem \ref{thm:matroid-intersection} also follows from the formula for the rank function of a union of two matroids, see~\cite[Theorem 11.3.1]{oxley}. The proof of the formula also uses Rado's theorem.
\end{remark}

\section{Amalgams of algebraic matroids over one-dimensional groups}\label{sec:amalgams}

In~\cite{BCD} and~\cite{evans1991projective} the authors use the notion of \emph{algebraically equivalent} realizations of an algebraic matroid.
Given some algebraic matroid $M(f_1,\dots,f_n)$ represented by $f_1,\dots,f_n \in L$. If there are $f_1',\dots,f_n' \in L$ such that
$\cl(f_i) = \cl(f_i')$ for all $i \in [n]$ then the matroid $M(f_1',\dots,f_n')$ equals $M(f_1,\dots,f_n)$ and we call the representations \emph{algebraically equivalent}. We write $f \sim f'$ for $f,f' \in L$ if $\cl(f) = \cl(f')$. The algebraic matroid represented by $f_1,\dots,f_n \in L$ is algebraically equivalent to a representation over a one-dimensional algebraic group $G$ if there are points $q_1,\dots,q_n \in G(L)$ that are the image of a generic point $p \in G(L)^d$ under a homomorphism of algebraic groups $\psi: G^d \to G^n$ such that $q_i \sim f_i$ for all $i \in [n]$.

Given a matroid $M_1$ realized by $(x,y) = (x_1,\dots,x_d,y_1,\dots,y_n)$ that is algebraically equivalent to a realization over $\mathbb{G}_a$ and a matroid $M_2$ realized by $(x,z) = (x_1,\dots,x_d,z_1,\dots,z_m)$ that is algebraically equivalent to a realization over a group $G$ which is either $\mathbb{G}_m$ or an \emph{elliptic curve}. Assume that both matroids have rank $d$ and the basis $x_1,\dots,x_d$. Let $M = M(x,y,z)$ be the algebraic matroid realized by $(x,y,z)$ with the ground set $E$. Let $\mathbb{E}_1$ and $\mathbb{E}_2$ be the endomorphism rings of $\mathbb{G}_a$ and $G$, with division rings $Q_1$ and $Q_2$. Then there are generic points $u\in\mathbb{G}_a(L)^d$ and $p\in G(L)^d$ and matrices $\alpha\in\mathbb{E}_1^{n\times d}$, $\beta\in\mathbb{E}_2^{m\times d}$ such that
\[u_i\sim x_i\sim p_i,\qquad y_j\sim\alpha_j(u),\qquad z_k\sim\beta_k(p).\]
Indeed, if $(x,y)\sim\psi(g)$, set $u_i=\psi_{x_i}(g)$ and $\psi_x$ for the matrix with rows $\psi_{x_1},\dots,\psi_{x_d}$; each row $\psi_{y_j}$ lies in the $Q_1$-span of the rows $\psi_{x_i}$, so by the Ore condition $e_j\psi_{y_j}=\alpha_j\psi_x$ for some $0\neq e_j\in\mathbb{E}_1$ and $\alpha_j\in\mathbb{E}_1^{d}$, and Lemma~\ref{lemm:matroids-over-groups}\ref{lemm:matroids-over-groups:1} gives $y_j\sim\alpha_j(u)$. The same argument applies to $(x,z)$. Set $q_i = (u_i,p_i) \in (\mathbb{G}_a \times G)(L)$ and denote with $C_i = Z(q_i) \subseteq \mathbb{G}_a \times G$ its \emph{locus}. With $\trdeg K(q_1,\dots,q_d) = d = \sum_i \trdeg K(q_i)$ it follows
\[ W \coloneqq Z((q_1,\dots,q_d)) = C_1 \times \dots \times C_d \subseteq (\mathbb{G}_a \times G)^d.\]
Fix a subset $S \subset E$ with $|S| = d$ and collect the $u_i$ and $\alpha_i$ appearing in $S$ into an $a \times d$ matrix $\alpha_S$ over $\mathbb{E}_1$. Collect all the $\beta_i$ appearing in $S$ into a $b \times d$ matrix $\beta_S$ over $\mathbb{E}_2$.

We define the following \emph{homomorphism of algebraic groups}
\[
\Phi_S\colon(\Ga\times G)^d\to\Ga^a\times G^b,\qquad \big((v_i,w_i)\big)_{i}\mapsto\big(\alpha_S \,v,\ \beta_S \, w\big).
\]

If we evaluate this morphism at the point $q = (q_1,\dots,q_d)$, we get with Lemma \ref{lemm:loci}
\[r_M(S) = \trdeg K(\Phi_S(q)) = \dim Z(\Phi_S(q)) = \dim \overline{\Phi_S(W)}\]
We set $X_j = \overline{\Phi_S(0 \times \dots  \times C_j \times \dots \times 0)}$ and since $\Phi_S$ is a morphism of algebraic groups we have
\[\overline{X_1 \oplus \dots \oplus X_d} = \overline{\Phi_S(W)}\]

For $I \subseteq [d]$ we set
\[X_I = \overline{\bigoplus_{i \in I} X_i}.\]

The following lemma follows directly from the theory of matroids over one-dimensional algebraic groups.
\begin{lemm}
\label{lemm:dim-of-matrix-image}
     Let $G$ be an \emph{elliptic curve} or $\mathbb{G}_m$. For $\alpha \in \mathbb{E}_1^{a \times d}$ we have $\dim \alpha(\mathbb{G}_a^d) = \rk_{Q_1} \alpha$ and for $\beta \in \mathbb{E}_2^{b \times d}$ we have $\dim \beta(G^d) = \rk_{Q_2} \beta$.
\end{lemm}

\begin{proposition}
    We have that $\dim \overline{\Phi_S(W)} = d$ if and only if for every $I \subseteq [d]$
    \[\rk_{Q_1} \alpha_{S,I} + \rk_{Q_2} \beta_{S,I} \geq |I|.\]
\end{proposition}
\begin{proof}
    Suppose the inequality does fail for some $I$. $X_I$ lies in the subgroup $\alpha_{S,I}(\Ga^I)\times\beta_{S,I}(G^I)$.
    And with Lemma \ref{lemm:dim-of-matrix-image} this subgroup has dimension $\rk_{Q_1} \alpha_{S,I} + \rk_{Q_2} \beta_{S,I} < |I|$. It follows
    \[\dim \overline{\Phi_S(W)} \leq \dim X_I + \dim X_{I^c} < |I| + |I^c| = d.\]
    Suppose $\dim \overline{\Phi_S(W)} < d$. We choose a minimal $I$ with $\dim X_I < |I|$. Hence we have $\dim {X_{I \setminus j}} = |I| -1$ for every $j \in I$ and $\dim X_I = |I| -1$.
    Fix some $j \in I$. We want to show that $X_j$ lies in a coset of a connected subgroup. Set $A \coloneqq {X_{I \setminus j}}$ and $Y \coloneqq X_j$. 
    $A$ is closed and irreducible as the image of an irreducible variety. We can assume that $Y$ is an irreducible curve, if $Y$ is only a point there is nothing to show. Now Lemma \ref{lemm:matroids-over-groups}\ref{lemm:matroids-over-groups:2} gives us $X_j \ominus X_j \subset \Stab (X_I)$. We fix some $x_j \in X_j(K)$. With $X_j \ominus x_j$ containing $0$ and being irreducible we have $0 \in X_j \ominus x_j \subseteq \Stab (X_I)^{\circ}$ which denotes the identity component of the stabilizer. Hence $X_j \subset x_j \oplus \Stab (X_I)^{\circ}$.
    We can do this now for every $j \in I$ and get $X_I \subseteq c \oplus \Stab (X_I)^\circ$ for $c = \bigoplus_{j \in I} x_j$, since $\Stab (X_I)^\circ$ is closed. But we also have that for any $c' \in X_I(K)$ that $c' \oplus \Stab (X_I)^\circ \subseteq X_I$. Thus $X_I = c \oplus \Stab (X_I)^\circ$. It follows $\dim \Stab (X_I)^\circ = |I| -1$.

    With Lemma \ref{lemm:matroids-over-groups}\ref{lemm:matroids-over-groups:3} follows $\Stab (X_I)^\circ = U \times V$ for closed subgroups $U \subseteq \mathbb{G}_a^a$ and $V \subseteq G^b$. We consider the projections $\pi_1: \mathbb{G}_a^a \times G^b \to \mathbb{G}_a^a$ and $\pi_2: \mathbb{G}_a^a \times G^b \to G^b$. Set
    \[A_I \coloneqq \alpha_{S,I}(\mathbb{G}_a^I) \subseteq \mathbb{G}_a^a, \  B_I \coloneqq \beta_{S,I}(G^I) \subseteq G^b.\]
    We have
    \[\overline{\pi_1(X_I)} = \overline{\alpha_{S,I}(\prod_{i\in I} \pi_1(C_i))},\]
    and as each $C_i$ maps dominantly into $\mathbb{G}_a$ we have that $\pi_1(X_I)$ is dense in $A_I$.
    We have that $X_I$ is a coset of $\Stab (X_I)^\circ$, hence $\pi_1(X_I) = \pi_1(c) + U$, which is already closed. Now
    \[\pi_1(c) + U = {\pi_1(X_I)} = A_I.\]
    As both $A_I$ and $U$ are closed connected subgroups we have $U = A_I$. The same argument gives us $V = B_I$ and it follows
    \[|I| -1 = \dim \Stab (X_I)^\circ = \dim A_I + \dim B_I = \rk_{Q_1} \alpha_{S,I} + \rk_{Q_2} \beta_{S,I}\]
    and the inequality fails for $I$.
\end{proof}

\begin{corollary}
\label{cor:intersection}
    If we use the notation of Section 3, then it follows with Theorem \ref{thm:matroid-intersection} that $M = M_1' \wedge M_2'$. 
\end{corollary}

\section{Algebraic equivalence of realizations}
We now want to show that the algebraic representations of certain matroids are all algebraically equivalent to a representation over a one-dimensional algebraic group $G$. This is already known for the \emph{Fano} and \emph{non-Fano} matroids and also for every \emph{projective geometry} in $L/K$, see~\cite[Theorem 3.3.1]{evans1991projective} and~\cite[Example 6.2]{BCD}. We will show the same for analogs of the \emph{dual Fano} and \emph{non-Fano} matroids over fields of characteristic $p$. For this the following lemmas will be helpful.

\subsection{Intersecting transcendental field extensions}
The statements are well known, but we provide short proofs.

\begin{lemm}
\label{lemm:rel-closure}
    Let $A \subseteq F$ be fields with $A$ algebraically closed and $t$ transcendental over $F$. Then $A(t)$ is algebraically closed in $F(t)$.
\end{lemm}
\begin{proof}
    We have that $A(t)$ and $F$ are free over $A$. With \cite[VIII. Theorem 4.3]{lang2012algebra}  and the fact that $F/A$ is regular the statement follows. 
\end{proof}

For subfields $A,B \subseteq L$ let $AB$ be the composite of $A$ and $B$ in $L$, that is the smallest subfield of $L$ which contains both $A$ and $B$. 

\begin{lemm}
\label{lemm:transc-intersection}
    Let $A,B \subseteq L$ be subfields. If $t \in L$ is transcendental over $AB$ then
    \[A(t) \cap B(t) = (A \cap B)(t).\]
\end{lemm}
\begin{proof}
    We only prove the inclusion $A(t) \cap B(t) \subseteq (A \cap B)(t)$. The other inclusion is trivial. Set $F = AB$. Every $f \in F(t)$ can be written uniquely as $f = \frac{p}{q}$ with $p, q \in F[t]$, $q$ monic, and $\gcd(p, q) = 1$ in $F[t]$.
    Now let $f \in A(t)$. Write $f = \frac{p}{q}$ with $p, q \in A[t]$, $q$ monic, and $\gcd(p, q) = 1$ in $A[t]$. Then $\gcd(p, q) = 1$ in $F[t]$ as well, and $\frac{p}{q}$ is the reduced form of $f$ over $F$.
    Hence the coefficients of the reduced form of $f$ lie in $A$. By the same argument, if $f \in B(t)$, they lie in $B$. Uniqueness of the reduced form means that for $f \in A(t) \cap B(t)$ these are the same coefficients. So they lie in $A \cap B$, and $f \in (A \cap B)(t)$.
\end{proof}

\begin{lemm}
\label{lemm:base-change}
    Let $A,B \subseteq L$ be algebraically closed subfields containing $K$. If $F = \cl(AB)$ and $t \in L$ is transcendental over $F$, then $\cl(A(t)) \cap \cl(B(t)) = \cl((A\cap B)(t))$.
\end{lemm}
\begin{proof}
    Choose some $z \in \cl(A(t)) \cap \cl(B(t))$ and let $f$ be its minimal polynomial over $F(t)$ and $g$ its minimal polynomial over $A(t)$. It follows that $f$ divides $g$ as polynomials over $F(t)$. Therefore all the coefficients of $f$ lie in $\cl(A(t)) \cap F(t)$. With Lemma \ref{lemm:rel-closure} we have $\cl(A(t)) \cap F(t) = A(t)$ and it follows that all the coefficients of $f$ lie already in $A(t)$.
    Similarly, we can get that $f \in B(t)[X]$. With Lemma \ref{lemm:transc-intersection} we have $A(t) \cap B(t) = (A \cap B)(t)$ and it follows that $z$ is algebraic over $(A\cap B)(t)$. Hence $z \in \cl((A \cap B)(t))$. The other inclusion is trivial.
\end{proof}

\subsection{The matroid \texorpdfstring{$\bm{M(K_4)}$}{M(K4)}}

Consider the matroid of the \emph{complete graph on $4$ vertices} $M(K_4)$. That is a matroid of rank $3$ with ground set $\{a,b,c,d,e,f\}$ and circuits $\{a,b,c\}$, $\{a,d,e\}$, $\{b,e,f\}$ and $\{c,d,f\}$. The matroid has the partition $\{a,f\},\{b,d\},\{c,e\}$ into points not on a \emph{common $3$-point line}. We will call a subset of the partition a \emph{matching}. With \cite[Theorem 2.1.2]{evans1991projective} follows that every algebraic realization of $M(K_4)$ is \emph{algebraically equivalent} to a realization over a one-dimensional algebraic group $G$. Using this one can show:
\begin{lemm}[\cite{evans1991projective} Corollary 3.2.1]
    Suppose $M(K_4)$ is realized by $a,b,c,d,e,f \in L$ as an algebraic matroid over an algebraically closed subfield $K_0 \subseteq L$. Then the lines $\cl(a,f)$ and $\cl(c,e)$ meet in a point over $K_0$. In other words, two matchings of $M(K_4)$ always meet in a point over $K_0$.
\end{lemm}

\begin{lemm}
\label{lemm:recog-MK4}
    Let $M$ be a \emph{simple matroid} of rank $3$ on $E = \{a,b,c,d,e,f\}$ where each circuit of $M(K_4)$ is a \emph{dependent} set in $M$. Call them $l_1,l_2,l_3,l_4$.
    Then those are the only dependent triples of $M$ and $M \cong M(K_4)$.
\end{lemm}
\begin{proof}
    This is well known, by the enumeration of all matroids of rank $3$ on $6$ elements, but we state a short proof.
    We notice that any three of the sets $l_1,l_2,l_3,l_4$ cover $E$.
    Let $F$ be a rank-$2$ flat of $M$ with $|F|\ge 3$. So $F$ contains two elements of some $l_i$. Since $M$ is simple and $l_i$ is dependent we have $l_i \subseteq F$.

    Suppose $F\ne l_i$, and pick $w\in F\setminus l_i$. $w$ is matched to exactly one element of $l_i$. So there are $s\ne t$ in $l_i$ such that $\{w,s\}$ lies in some different $l_j$ and $\{w,t\}$ lies in some different $l_k$. As before we have $l_j,l_k \subseteq F$. The three sets $l_i,l_j,l_k$ are distinct, so $E\subseteq F$. This contradicts $r(M)=3$. Hence $F=l_i$.
\end{proof}

\subsection{The dual Lazarson matroids are only realizable over algebraic groups}\label{sec:alg-reps}

If we extend the field $K$ by some elements $a_1,\dots,a_n$ we will always assume that $K(a_1,\dots,a_n)$ is already algebraically closed from now on in this section.

Given variables $x_0,\dots, x_{p+1}$. For $i \in [p+1]$ set $y_i = \sum_{j \neq i} x_j$. Let $L_p^*, \NL_p^*$ be the matroids realized by $\{x_0,\dots,x_{p+1},y_1,\dots,y_{p+1}\}$ over $\mathbb{F}_{p}$ or $\QQ$. The matroid $L_p^*$ is the dual of the so-called \emph{Lazarson matroid} $L_p$ and the matroid $\NL_p^*$ is the dual of the \emph{non-Lazarson matroid} $\NL_p$. The two matroids can be seen as generalizations of the \emph{Fano} and \emph{non-Fano} matroids to characteristic $p$. $L_p$ is only algebraic over fields of characteristic $p$ and $\NL_p$ is only linear over fields of characteristic not equal to $p$. 

\begin{theorem}
\label{thm:alg-equivalence}
    Let $N_p$ be either $L_p^*$ or $\NL_p^*$ and realized by some arbitrary elements
    \[x_0,\dots,x_{p+1},y_1,\dots,y_{p+1} \in L\] 
    as an algebraic matroid over $K$.
    Then there is an \emph{algebraically equivalent} realization over some \emph{connected one-dimensional algebraic group} $H$.
    Furthermore for $N_p = L_p^*$ we have $H \cong \mathbb{G}_a$ and for $N_p = \NL_p^*$ we have $H \cong \mathbb{G}_m$ or $H$ is some elliptic curve.
\end{theorem}
\begin{proof}
    We will construct several realizations of $M(K_4)$ over $K$ from $N_p$ and glue them together into an algebraically equivalent realization of $N_p$ over a one-dimensional algebraic group.
    
    Fix $i,j \in [p+1]$. Set $S_{ij} = \{x_1,\dots, x_{p+1}\} \setminus \{x_i,x_j\}$. Then the \emph{contraction} $N_p/S_{ij}$ has $p+4$ elements and all the $y_k$ with $y_k \neq y_i$ and $y_k \neq y_j$ are \emph{parallel} to each other; we represent them by $z$. Now consider the \emph{simplification} $\{x_0,x_i,x_j,y_i,y_j,z\}$. We have the \emph{lines} $\{x_0,x_j,y_i\}$, $\{x_0,x_i,y_j\}$, $\{x_i,y_i,z\}$ and $\{x_j,y_j,z\}$. Hence this matroid is isomorphic to $M(K_4)$ with the \emph{matchings} $\{x_0,z\}$, $\{x_i,x_j\}$ and $\{y_i,y_j\}$. Therefore we get the point
    \[K(S_{ij})(p_{ij}) = K(S_{ij})(x_i,x_j) \cap K(S_{ij})(y_i,y_j)\]
    over $K(S_{ij})$. As $\{x_i,x_j,y_i,y_j\}$ is a flat in $N_p$ and $S_{ij}$ is independent over $K(x_i,x_j,y_i,y_j)$ we get with Lemma \ref{lemm:base-change}
    \[K(x_i,x_j)(S_{ij}) \cap K(y_i,y_j)(S_{ij}) = (K(x_i,x_j) \cap K(y_i,y_j))(S_{ij})\]
    has transcendence degree $p$ over $K$ and we can take $p_{ij}$ to be \emph{points} over $K$.
    As the intersection of the flats spanned by $\{x_i,x_j\}$ and $\{y_i,y_j\}$ is empty in $N_p$ the elements $p_{ij}$ are not parallel to any element of $N_p$.
    Assume that $p_{ij}$ is parallel to $p_{kl}$. It follows $p_{ij} \in K(x_i,x_j) \cap K(x_k,x_l)$. But this can only happen if $\{i,j\} = \{k,l\}$. Hence any two points are independent.
    
    Now set 
    \[\Pi_{ijk} = \{x_i,x_j,x_k,p_{ij},p_{ik},p_{jk}\}\]
    and
    \[\Pi'_{ijk} = \{y_i,y_j,y_k,p_{ij},p_{ik},p_{jk}\}.\]
    We have that $\{p_{ij},p_{ik},p_{jk}\} \subseteq K(x_i,x_j,x_k)$ and $\{p_{ij},p_{ik},p_{jk}\} \subseteq K(y_i,y_j,y_k)$. These two fields are not equal and both have transcendence degree $3$. Hence $\{p_{ij},p_{ik},p_{jk}\}$ has rank $2$ and is a line. In addition we have the lines
    \[\{x_i,x_j,p_{ij}\},\{x_i,x_k,p_{ik}\},\{x_j,x_k,p_{jk}\}\]
    and
    \[\{y_i,y_j,p_{ij}\},\{y_i,y_k,p_{ik}\},\{y_j,y_k,p_{jk}\}.\]
    It follows with Lemma \ref{lemm:recog-MK4} that $\Pi_{ijk}$ and $\Pi'_{ijk}$ both realize $M(K_4)$ over $K$. By~\cite[Theorem 2.1.2]{evans1991projective} each of them is algebraically equivalent to a realization over a connected one-dimensional algebraic group, and by projective uniqueness of $M(K_4)$ over division rings we can fix the realization: for $\Pi_{ijk}$ there are a group $G$ and algebraically independent points $s_i,s_j,s_k\in G(L)$ with
    \[x_l\sim s_l\ (l\in\{i,j,k\}),\qquad p_{lm}\sim s_l\ominus s_m\ (l,m\in\{i,j,k\}),\]
    and similarly for $\Pi'_{ijk}$ with $y_l$ in place of $x_l$.

    We repeatedly use the following observation. Let $\pi\colon H'\to H$ be an isogeny, $a\in H(L)$, $c\in H(K)$ and $a'\in H'(L)$ with $\pi(a')=a\oplus c$. Then $a'\sim a$, since $\pi$ is finite and $K(a\oplus c)=K(a)$. Such lifts $a'$ exist because $\pi$ is surjective and $L$ is algebraically closed.

    \emph{Step 1: the elements $x_1,\dots,x_{p+1}$ and $p_{ij}$.} We show by induction on $k\in[3,p+1]$ that there are a connected one-dimensional algebraic group $H_k$ and points $p_1,\dots,p_k\in H_k(L)$ with
    \[x_i\sim p_i \text{ for } i\in[k],\qquad p_{ij}\sim p_i\ominus p_j \text{ for } i\in\{1,2\},\ j\in[k]\setminus\{i\}.\]
    For $k=3$ this is the realization of $\Pi_{123}$. For the induction step, take the realization $q_1,q_2,q_{k+1}\in G(L)$ of $\Pi_{1,2,k+1}$, with $x_{k+1}\sim q_{k+1}$. Apply Proposition~\ref{lemm-gauge} to $a=(p_1,q_1)$ and $b=(\ominus p_2,\ominus q_2)$: the points $p_1,p_2$ are independent, and $p_1\ominus p_2\sim p_{12}\sim q_1\ominus q_2$. We obtain a group $H_{k+1}\subseteq H_k\times G$ such that the projections $\pi\colon H_{k+1}\to H_k$ and $\pi'\colon H_{k+1}\to G$ are isogenies, together with $r_1,r_2\in H_{k+1}(L)$ and constants $c_1,c_2$ with $(p_1,q_1)=r_1\oplus c_1$ and $( p_2, q_2)= r_2\ominus c_2$. Choose $r_i\in H_{k+1}(L)$ with $\pi(r_i)=p_i$ for $3\le i\le k$, and $r_{k+1}$ with $\pi'(r_{k+1})=q_{k+1}$. Since $\pi(r_i)$ and $\pi'(r_i)$ differ from $p_i$ and $q_i$ only by constants, the observation gives $r_i\sim x_i$ for all $i\in[k+1]$, $r_i\ominus r_j\sim p_i\ominus p_j\sim p_{ij}$ for $i\in\{1,2\}$, $j\le k$, and $r_i\ominus r_{k+1}\sim q_i\ominus q_{k+1}\sim p_{i,k+1}$ for $i\in\{1,2\}$. Rename each $r_i$ to $p_i$.

    \emph{Step 2: the elements $y_1,y_2,y_3$.} Let $p_1,\dots,p_{p+1}\in H_{p+1}(L)$ be as in Step~1, and let $t_1,t_2,t_3\in G'(L)$ be the realization of $\Pi'_{123}$. Apply Proposition~\ref{lemm-gauge} to $a=(p_1\ominus p_2,\,t_1\ominus t_2)$ and $b=(p_2\ominus p_3,\,t_2\ominus t_3)$: these realize $p_{12}$ and $p_{23}$, which are independent, and $a\oplus b$ realizes $p_{13}$ in both factors. We obtain $H\subseteq H_{p+1}\times G'$ with isogenies $\pi\colon H\to H_{p+1}$, $\pi'\colon H\to G'$ and points $p_a,p_b\in H(L)$ that project to the components of $a$ and $b$ up to constants. Choose $p'_2$ with $\pi(p'_2)=p_2$, set $p'_1=p_a\oplus p'_2$ and $p'_3=p'_2\ominus p_b$, choose $p'_i$ with $\pi(p'_i)=p_i$ for $i\ge4$, and choose $t_1'\in H(L)$ with $\pi'(t'_1)=t_1$. Then $\pi(p'_i)=p_i$ up to constants for all $i$, so the relations of Step~1 still hold for the $p'_i$. Moreover, up to constants,
    \[\pi'(t'_1\ominus p'_1\oplus p'_2)=t_1\ominus\pi'(p_a)=t_2,\qquad \pi'(t_1'\ominus p'_1\oplus p'_3)=t_1\ominus\pi'(p_a)\ominus\pi'(p_b)=t_3,\]
    so $y_1\sim t_1'$, $y_2\sim t_1'\ominus p'_1\oplus p'_2$ and $y_3\sim t_1'\ominus p'_1\oplus p'_3$.

    \emph{Step 3: the remaining elements.} Set $h=(p'_1,\dots,p'_{p+1},t_1')\in H(L)^{p+2}$. Since $\{x_1,\dots,x_{p+1},y_1\}$ is a basis of $N_p$, the point $h$ is generic. So far we have found $v\in\mathbb{E}^{p+2}$ with entries in $\{0,\pm1\}$ such that $f\sim v(h)$ for every $f\in\{x_1,\dots,x_{p+1},y_1,y_2,y_3\}\cup\{p_{1k},p_{2k}\}$. For $k \ge 4$, the element $y_k$ lies on the two distinct lines $\{y_1,p_{1k},y_k\}$ and $\{y_2,p_{2k},y_k\}$, so $\cl(y_k)=\cl(y_1,p_{1k})\cap\cl(y_2,p_{2k})$. By Lemma~\ref{lemm-modular-intersection} this intersection is $\cl(w(h))$ for some $w\in\mathbb{E}^{p+2}$, hence $y_k\sim w(h)$. Similarly, with $F_j=\{y_j\}\cup\{x_i\mid i\neq 0,j\}$ we have $\cl(x_0)=\bigcap_{j}\cl(F_j)$, and applying Lemma~\ref{lemm-modular-intersection} inductively gives $x_0\sim w_0(h)$ for some $w_0\in\mathbb{E}^{p+2}$. Collecting the rows $v$, $w$ and $w_0$ into a matrix $\psi\in\mathbb{E}^{(2p+3)\times(p+2)}$ gives a homomorphism $\psi\colon H^{p+2}\to H^{2p+3}$ with $\psi(h)\sim(x_0,\dots,x_{p+1},y_1,\dots,y_{p+1})$ coordinatewise. Hence $N_p$ is algebraically equivalent to a realization over $H$.
    In particular, for $N_p = L_p^*$, $H$ needs to be $\mathbb{G}_a$ as $L_p^*$ is only linear over division rings of characteristic $p$. For $N_p = \NL_p^*$, $H$ needs to be $\mathbb{G}_m$ or an elliptic curve as $\NL_p^*$ is only linear over division rings of characteristic not $p$.
\end{proof}

\section{Algebraic matroids which are not coalgebraic}

The \emph{Lazarson matroid} $L_p$ over $\overline{\FF}_p$ is realized by
\[x_1,\dots,x_{p+1},y_1,\dots,y_{p+1},z\]
where $z = \sum_{i=1}^{p+1} x_i$ and $y_i = z - x_i$.
The \emph{non-Lazarson matroid} $\NL_p$ over $\overline{\FF}_p$ is realized by
\[u_1,\dots,u_{p+1},v_1,\dots,v_{p+1},w\]
where $w = \prod_{i=1}^{p+1} u_i$ and $v_i = \frac{w}{u_i}$. Here the elements $u_1,\dots,u_{p+1},x_1,\dots,x_{p+1}$ are algebraically independent over $\overline{\FF}_p$.
That these matroids are actually \emph{dual} to $L_p^*$ and $\NL_p^*$ is a routine check.
Consider the \emph{union} $M$ of the following two matroids over $\overline{\FF}_p$:
\[\tilde{L}_p = M(x_1,\dots,x_{p+1},0,\dots,0,y_1,\dots,y_{p},y_{p+1},z)\]
and
\[\widetilde{\NL}_p = M(0,\dots,0,u_1,\dots,u_{p+1},v_1,\dots,v_{p},w,v_{p+1})\]
with the ground set $\{a_1,\dots,a_{p+1},b_1,\dots,b_{p+1},c_{1},\dots,c_{p+2}\}$.
Then \\ $\{c_1,\dots,c_p,c_{p+2}\} \cup \{b_1,\dots,b_p,c_{p+1}\}$ is a basis of $M$, as both are disjoint bases of $\tilde{L}_p$ and $\widetilde{\NL}_p$. 

The algebraic matroids in the following theorem are inspired by the counterexample of Larson and Le in \cite{larson2026duals}.
\begin{theorem}
    Consider the matroid $M_p$ represented by
    \[x_1,\dots,x_{p+1},u_1,\dots,u_{p+1},v_1y_1,\dots,v_p y_{p}, w \cdot y_{p+1}, v_{p+1} \cdot z\]
    over $\overline{\FF}_p$. Then the dual of $M_p$ is not algebraic over any field.
\end{theorem}
\begin{proof}
    The subset $\{b_1,\dots,b_p, c_1,\dots,c_p, c_{p+1}\}$ is already dependent. We can see this by using $u_1,\dots,u_p$ to transform 
    \[v_1y_1,\dots,v_p y_{p}, w y_{p+1}\]
    into
    \[u_{p+1}y_1,\dots,u_{p+1} y_{p}, u_{p+1}  y_{p+1}.\]
    Now we have
    \[u_{p+1} \sum_{i=1}^{p+1} y_i =u_{p+1} \sum_{i=1}^{p+1} (z - x_i) = u_{p+1} (z - \sum_{i=1}^{p+1} x_i) = 0\]
    and the set $\{b_1,\dots,b_p, c_1,\dots,c_p, c_{p+1}\}$ is dependent in $M_p$. If we contract the $x_i$ in $M_p$ the matroid $M_p/x$ is equal to $\widetilde{\NL}_p$ and if we contract the $u_i$ the matroid $M_p/u$ is equal to $\tilde{L}_p$. Therefore $M_p$ is not the union of $M_p/u$ and $M_p/x$ where we leave the contracted elements as loops. Assume the dual matroid $M_p^*$ is algebraic. Then $M_p^*$ is algebraic over a field $K$ of characteristic $p$, as the dual of the Lazarson matroid is only algebraic over a field of characteristic $p$. Then the $c_i$ form a basis of $M_p^*$, $M_p^*\setminus a = \widetilde{\NL}_p^*$ and $M_p^*\setminus b = \tilde{L}_p^*$. With Theorem \ref{thm:alg-equivalence} every algebraic realization of $M_p^*\setminus a$ is algebraically equivalent to a realization over $\mathbb{G}_m$ or an elliptic curve and every algebraic realization of $M_p^* \setminus b$ is algebraically equivalent to a realization over $\mathbb{G}_a$. With Corollary \ref{cor:intersection} it follows that $M_p^*$ must be the intersection of $M_p^*\setminus a$ and $M_p^*\setminus b$ where we leave the deleted elements as coloops. This contradicts that $M_p$ is not the union of $M_p/u$ and $M_p/x$ where we leave the contracted elements as loops. It follows that the dual of $M_p$ is not algebraic over any field and proves Theorem \ref{thm:main}. 
\end{proof}
Our main theorem now follows directly.
\begin{theorem}
\label{thm:counterexample}
	Over any field $K$ of \emph{positive characteristic} there is an \emph{algebraic matroid} $M$ over $K$ such that the \emph{dual matroid} $M^*$ is \emph{not algebraic} over any field.
\end{theorem}

\begin{remark}
    There are several matroids similar to $M_p$ with non-algebraic dual. For example one could use different twists to glue the matroids $L_p^*$ and $\NL_{p}^*$ together. The twist we used just seems to be the easiest one that works. It does not work to glue $L_p$ and $\NL_p$ together with some twist and arrive at the same result. 
\end{remark}

\section{Volume polynomials}\label{sec:volume}

In \cite{castillo2020multidegrees} the authors define algebraic polymatroids. 
Suppose we have a chain of field extensions $K \subseteq L_i \subseteq L$ for all $i\in E$ for some finite set $E$.
Denote with $L_A$ the composite of field extensions $L_i \subseteq L$ for $i \in A$. This is the smallest subfield of $L$ containing $L_i$ for all $i \in A$. 

\begin{definition}
    Let $\mathcal{P} = (E,r)$ be an integral polymatroid on the ground set $E$ with rank function $r$.
    We call $\mathcal{P}$ \emph{algebraic over $K$} if there are field extensions $K \subseteq L_i \subseteq L$ for $i \in E$ such that
    \[r(A) = \trdeg_K(L_A) \text{ for all } A \subseteq E.\]
\end{definition}

See~\cite[Chapter 4]{murota} for several equivalent definitions of integral polymatroids.
We view integral polymatroids and $M$-\emph{convex} sets as the same objects and denote with $\mathcal{P}^*$ the dual of a polymatroid $\mathcal{P}$.

Algebraic polymatroids are connected to a special family of polynomials coming from projective geometry called \emph{volume polynomials}.
These are a subclass of the so-called \emph{Lorentzian polynomials}, see \cite{branden2020lorentzian}.
A Lorentzian polynomial is a homogeneous polynomial with non-negative coefficients satisfying a certain log-concavity property on the positive orthant. Given a projective variety $X$ over a field $K$ of dimension $d$ and a collection $D = (D_1,\dots,D_n)$ of semiample divisors. The \emph{realizable volume polynomial} of $D$ is defined as
\[f_D(x_1,\dots,x_n) = \frac{1}{d!} \int_X \left( \sum_{i=1}^n x_iD_i \right)^d.\]
This is a homogeneous polynomial of degree $d$. For a survey on volume polynomials, see \cite{huh2026volumepolynomials}. The coefficients of volume polynomials obey additional inequalities compared to Lorentzian polynomials. While the support of a Lorentzian polynomial can be any polymatroid, we have the following restriction for volume polynomials:

\begin{proposition}[\cite{GHMSW2025preservingvolume}, Proposition 5.4]
    The following statements are equivalent for any polymatroid $\mathcal{P}$ with set of bases $J$ and any field $K$.
    \begin{enumerate}
        \item There is a realizable volume polynomial over $K$ whose support is $J$.
        \item The polymatroid $\mathcal{P}$ is \emph{algebraic over $K$}.
    \end{enumerate}
\end{proposition}

Theorem~\ref{thm:main} immediately implies that the class of algebraic polymatroids is also not closed under duality. 

Matroid duality can be generalized to the class of volume polynomials.
Consider the dual pair of polynomial rings $ \RR[x] = \RR[x_1,\dots,x_n]$ and $\RR[\partial] = \RR[\partial_1,\dots,\partial_n]$ where $\RR[\partial]$ acts on $\RR[x]$ via the usual rule of partial differentiation. Denote with $\RR[\partial]_{\leq \mu}$ all polynomials in $\RR[\partial]$ of multidegree at most $\mu \in \ZZ^n$. Then $g \in \RR[\partial]_{\leq \mu}$ is a \emph{realizable covolume polynomial} over $K$ if and only if $g(\partial) \circ \left( \frac{1}{\mu!}x^\mu \right)$ is a realizable volume polynomial over $K$. This definition does not depend on $\mu$.
The above correspondence then yields an analogous statement for volume polynomials and their duals:

\begin{theorem}
\label{thm:volume}
    Let $K$ be a field of positive characteristic. There exists a realizable covolume polynomial $g$ over $K$ such that the support of $g$ is not the support of any realizable volume polynomial.
\end{theorem}
\begin{proof}
    Let $g\in \RR[\partial]_{\leq \mu}$ be a realizable covolume polynomial such that the support of the dual $f = g(\partial) \circ \left( \frac{1}{\mu!}x^\mu \right)$ is the set of bases of an algebraic polymatroid $\mathcal{P}$ with non-algebraic dual $\mathcal{P}^*$. The dualization of $\mathcal{P}^*$ to $\mathcal{P}$ with respect to $\mu$ is the same operation as taking the support of the dual $f$ of $g$ with respect to $\mu$. Hence, the support of $g$ is not the support of any realizable volume polynomial.
\end{proof}

\printbibliography
\end{document}